\documentclass[11pt]{amsart}
\usepackage{amsmath}
\usepackage{amssymb}
\usepackage{amsthm}
 \DeclareMathOperator{\Zeros}{Zeros}
 \DeclareMathOperator{\Ideal}{Ideal}
\DeclareMathOperator{\Bl}{B\l} \DeclareMathOperator{\bl}{b\l}
\def\inv{^{-1}}

\newcommand{\bsl}{\begin{slide}{}}

\newcommand{\A}{\mathcal{A}}

\newcommand{\X}{\mathcal{X}}

\newcommand{\Cal}[1]{\mathcal #1}

\def\sbr #1.{^{[#1]}}
\def\sfl #1.{^{( #1)}}
\def\spr #1.{^{( #1)}}

\def\inv{^{-1}}
\def\?{{\bf{??}}}

\def\Hilb{\text{Hilb}}

\def\Proj{\text{Proj}}

\def\A{\Bbb A}
\def\X{\Cal X}

\def\C{\Bbb C}

\def\P{\mathbb P}

\def\ord{\text{\rm ord}}

\def\Spec{\text{\rm Spec} }

\def\ss{\vskip.12in}
\def\ms{\vskip.1in}
\def\Q{\Bbb Q}

\def\O{\Cal O}

\def\y{\bar{y}}

\def\Sym{\text{Sym}}

\def\rk{\text{rk}}

\def\1/2{\frac{1}{2}}

\def\I{\Cal I}

\def\2{{[2]}}
\def\l{\ell}
\def\nl{\newline}
\newtheorem*{thm*}{Theorem}
\newtheorem*{mainthm}{Main Theorem}
\newtheorem{thm}{Theorem}[section]
\newtheorem{cor}{Corollary}[section]

\newtheorem{prop}[cor]{Proposition}
\newtheorem{rem}[cor]{Remark}
\begin{document}
\title{Cycle map on Hilbert Schemes of Nodal Curves }
\author{Ziv Ran}
\address{ Mathematics Department, UC Riverside, CA 92521, USA}
\email{ziv.ran@ucr.edu}
\thanks{Research supported in part by the NSA under grant MDA 904-02-1-0094}
 \abstract We
study the structure of the relative Hilbert scheme for a family of
nodal (or smooth) curves via its natural cycle map  to the relative
symmetric product. We show that the cycle map is the blowing up of
the discriminant locus, which consists of cycles with multiple
points. We discuss some applications and connections, notably with
birational geometry and intersection theory on Hilbert schemes of
smooth surfaces.
\endabstract
\maketitle

\section*{Introduction} An object of  central importance in
classical algebraic geometry is a family of projective curves,
given by a projective morphism $$\pi:X\to B$$ with smooth general
fibre. One wants to take $B$ itself projective, which means one
must allow some singular fibres. We will assume our singular
fibres al all nodal. Of course, by semistable reduction, etc., any
family can be modified so as to have this property, without
changing the general fibre $X_b=\pi\inv(b)$. Many questions of
classical geometry involve point-configurations on fibres $X_b$
with $b\in B$ variable. From a modern standpoint, this means they
involve the relative Hilbert scheme
$$X\sbr m._B =\Hilb_m(X/B).$$ This motivates the interest in
studying $X\sbr m._B$ and setting it up as a tool for studying the
geometry, e.g. enumerative geometry, of the family $X/B$. This paper
is a step in this direction. Our focus will be on the \emph{cycle
map} (sometimes called the 'Hilb-to-Chow' map) $\frak c_m$, which in
this case takes values in the relative symmetric product $X\spr
m._B$.\par Our main result is that $\frak c_m$ is the blowing-up of
the discriminant locus in $X\spr m._B$. This result will be proven
in \S 1. As we shall see, the proof amounts to a fairly complete
study, locally over $X\spr m._B$, of $\frak c_m.$ We shall see in
particular that $\frak c_m$ is a small resolution of singularities;
in fact in 'most' cases the non-point fibres of $\frak c_m$ are
chains of rational curves (with at most $m-1$ components). In \S 2
we will consider applications of the result of \S 1 to the further
study of $X\sbr m._B$ and $\frak c_m.$ We will give a formula for
the canonical bundle of $X\sbr m._B$ showing that $\frak c_m$ 'looks
like' a flipping contraction; in fact, $\frak c_2$ is none other
than the \emph{Francia flip} and admits a natural 2:1 covering by
the flop associated to a 3-fold ODP
. We will also give a simple
formula for the \emph{Euler number} of $X\sbr m._B$. In \S3 we will
discuss the Chern classes of \emph{tautological bundles}. These are
bundles whose fibre at a point representing a scheme $z$ is
$H^0(E\otimes\O_z)$, where $E$ is a fixed vector bundle on $X$.\par
This paper has substantial intersection with the Author's papers
\cite{R,R2,R3} where some of the results are proven in greater
detail.\par As to the relevance of this paper to the theme of
'projective varieties with unexpected properties' I can only say
that the close links- some exposed below- of the Hilbert scheme, a
priori a purely algebraic object, to classical projective geometry
were quite unexpected by me, though this is probably due only to my
own ignorance.
\subsection*{Acknowledgments} A preliminary version of this work
was presented at the Siena conference in June '04. I would like to
thank the conference's organizers, especially Luca Chiantini, for
their hard and successful work putting together this memorable and
valuable mathematical event. I would also like to thank the
participants of both  the Siena conference and a subsequent one in
Hsinshu, Taiwan, especially Rahul Pandharipande and Lih-chung Wang,
for valuable input into \S3.
\section{The cycle map as blow-up}
Our main object of study is family of projective curves
$$\pi: X\to B$$
whose fibres $X_b=\pi\inv(b)$ are smooth for $b\in B$ general. We
shall make the following \begin{center} \emph {Essential
hypothesis: $X_b$ is
 nodal for all $b\in B.$}\end{center}\par  We shall also make
the (nonessential, but convenient) hypothesis that $X,B$ are smooth
of dimension 2,1 respectively.\par Geometry of the family largely
amounts to the study of families of subvarieties (more precisely
{\it{subschemes}})
$$\{Z_b\subset X_b, b\in B\}$$ of some fixed degree (length) $m$ over
$B$.

The canonical parameter space for subschemes is the {\emph{relative
Hilbert scheme}} $$X\sbr m._B =\Hilb_m(X/B).$$So (ordinary)points
$z\in X\sbr m._B$ correspond 1-1 with pairs $(b,Z)$ where $b\in B$
and $Z\subset X_b$ is a length-$m$ subscheme. More generally, for
any artin local $\C$-algebra $R$ and $S=\Spec(R)$, we have a
bijection between diagrams
$$\begin{matrix} S&\stackrel{f}{\rightarrow} & X\sbr
m._B\\&\stackrel{\searrow}{f_0}&\downarrow\\&&B\end{matrix}$$ and
$$\begin{matrix}
Z&\subset&X_S&\to&X\\&\searrow&\downarrow&&\downarrow\\
&&S&\stackrel{f_0}{\to}&B\end{matrix} $$ with the right square
cartesian and $Z/S$ flat of relative length $m$.

As usual in Algebraic Geometry, we study a complex object like
$X\sbr m._B$ by relating it (mapping it) to other (simpler ?)
objects. One approach (not pursued here, but see \cite{R}) is to
relate $X\sbr m._B$ (albeit only by correspondence, not morphism) to
$X\sbr m-1._B.$ This leads to studying {\emph{flag Hilbert
schemes}}. These have a rich  geometry; they are generally
singular.\nl We focus here on another approach, based on the
{\bf{cycle map}}
$$\frak c_m:X\sbr m._B\to \Sym^m_B(X)=:X\sfl m._B$$
$$Z\mapsto \sum\limits_{p\in X}{\text{length}}_p(Z)p.$$ Clearly,
$\frak c_m$ is an iso off the locus of cycles whose support meets
the critical or singular locus\begin{center}  sing$(\pi)=$ locus in
$X$ of singular points of fibres of $\pi.$\end{center}

 \begin{mainthm} $\frak c_m$ is the blowing-up of the
{\bf{discriminant locus}}  $$D^m=\{\sum m_ip_i : \exists
m_i>1\}\subset X\sfl m._B$$
\end{mainthm}

 Recall that if $I$ is an ideal on scheme $X$, we have a surjection
of graded algebras from the symmetric algebra on $I$ to the Rees or
blow-up algbera
$$\Sym^\bullet(I)\to \bigoplus\limits_0^\infty I^j$$ Applying the Proj
functor, we get a closed embedding (maybe strict) of schemes over
$X$
$$\Bl_I(X)\subseteq\P(I)$$ of the blow-up into the 'singular
projective bundle' $\P(I)$, whose fibres over $X$ are projective
spaces of varying dimensions. Note that $\P(I)$ may be reducible,
while $\Bl_I(X)$ is always an integral scheme if $X$ is.
Concretely, these schemes may be described, locally over $X$, as
follows: if $f_1,...,f_r$ generate $I$, take formal homogeneous
coordinates $T_1,...,T_r$, then as subschemes of
$X\times\P^{r-1}$,
\begin{eqnarray*}
\Bl_I(X)=\Zeros(
G(T_1,...,T_r):G(f_1,...,f_r)=0, G\ {\text{homogeneous}})~~~~~\\
\P(I)=\Zeros(
 G(T_1,...,T_r):G(f_1,...,f_r)=0, G\ {\text{homogeneous linear}})~~~~\\
\end{eqnarray*}
Thus, the inclusion $\Bl_I(X)\subseteq\P(I)$ is strict iff $I$
admits nonlinear syzygies; the case of the discriminant locus, to
be studied below, will provide examples of such ideals.
\ss\noindent\emph{Remark} Will see in the proof that\newline
$\bullet\ \ X\sbr m._B$ is smooth (over $\C$) of dimension $m+\dim
B$\newline$\bullet\ \ \frak c_m$ is a small map (in fact, if each
$X_b$ has at most $\nu$ nodes-- usually, $\nu=1\ $-- then fibres
of $\frak c_m$ have dimension at most $\min(\nu, m/2)$).\par
Clearly, $D^m$ is a prime Weil divisor on $X\sfl m._B$, in fact
$$D^m\sim_{bir} X\times_BX\sfl m-2._B$$ because a general $z\in D$
has the form $$z=2p_1+p_2+...+p_{m-1}$$.
On the other hand, near cycles meeting
sing$(\pi)$, esp. 'maximally singular' cycles
$$z=mp, p\in{\text{sing}}(\pi),$$ it's not clear a priori what (or how
many) defining equations $D^m$ has (the proof below will yield a
posteriori $m$ minimal equations locally at maximally singular
cycles).\nl Note that locally at maximally singular cycles, the
relative Cartesian product $X^m_B$ is a complete intersection with
equation $x_1y_1=...=x_my_m$, with the projection to $B$ given by
$t=x_1y_1$, while $X\sfl m._B$ is a quotient of a complete
intersection
$$(x_1y_1=...=x_my_m)/\text{ symmetric group}\ \frak S_m.$$
We will see that $X\spr m._B$
 is
not $\Q$-factorial: in fact,  $D^m$ is not $\Q$-Cartier;\nl Worse,
$X\sfl m._B$ is not even $\Q$-Gorenstein: we shall see that it
admits a small discrepant resolution $X\sbr m._B$.\nl Nonetheless,
being quotient by a finite group and smooth in codimension 1,
 $X\sfl m._B$ is normal and Cohen-Macaulay.
\ss The {\emph {plan of proof}} is as follows.\newline $\diamond\ $
Construct explicit (analytic) model of $X\sbr m._B$ and $\frak c_m$,
locally over $X\sfl m._B$; in particular, conclude that $X\sbr m._B$
is smooth and $\frak c_m$ is small, so $\frak c_m\inv(D^m)$ is
Cartier divisor.\nl
$\diamond\ $ The Universal property of blowing up now yields a
factorization
$$\begin{matrix}X\sbr m._B&\stackrel{\frak c'_m}{\to}&\Bl_{D^m}X\sfl
m._B\\&\stackrel{\searrow}{\frak c_m}&\downarrow \bl\\&&X\sfl
m._B\end{matrix}$$ Then we check locally (over the blowup) that
$\frak c'_m$ is an iso.\ss To start the proof, fix an analytic
neighborhood $U$ of fibre a node $p$, so the family is given in
local analytic coordinates by $xy=t.$\par

 For the local  study, the first question is:
what are fibres of $\frak c_m$?\newline Now locally in the \'etale
topology, all fibres are (essentially) products of fibres $\frak
c_{m_i}\inv(m_ip_i)$. So suffices to study
$$\frak c_m\inv(mp), p\in{\text{sing}}(\pi).$$ Then,
$$\frak c_m\inv(mp)=\Hilb_m^0(R)$$where $R$ is the formal
power series ring $$R=\C[[x,y]]/(xy).$$ Here $\Hilb_m^0$ denotes
the punctual Hilbert scheme.
\begin{prop} $\Hilb_m^0(R)$ is a chain
of $m-1$ smooth rational curves meeting normally
$$C^m_1\cup_{q^m_2}...\cup_{q^m_{m-1}} C^m_{m-1}:$$
$$\begin{matrix}&&Q^m_2&&Q^m_{m-1}&&
\\ \\&C^m_1&&C^m_2\hskip 1cm ...\hskip 1cm &&C^m_{m-1}&\\ \\
Q^m_1&&&&&&Q^m_m
\end{matrix}$$
\begin{center} Fig. 1\end{center}\vskip 1cm $q^m_i=(x^{m+1-i}, y^i),$
$$C^m_i\setminus\{q^m_i,
q^m_{i+1}\}=\{I^m_i(a)=(ax^{m-i}+y^i):a\neq 0\}$$\end{prop} NB
$\lim\limits_{a\to 0} I^m_i(a)=q^m_i,\lim\limits_{a\to \infty}
I^m_i(a)=q^m_{i+1}$.\proof See \cite{R2}\endproof\par Given this,
the next question is: what does the full Hilbert scheme look like
along $\Hilb^0$, e.g. locally near $q^m_i$?
\begin{prop} The universal flat deformation of the ideal
$q^m_i=(x^{m+1-i}, y^i), i=1,...,m$, rel $B$, is $(f,g)$ where
$$f=x^{m+1-i}+f^1_{m-i}(x)+vy^{i-1}+f^2_{i-2}(y),$$
$$g=y^i+g^1_{i-1}(y)+ux^{m-i}+g^2_{m-i-1}(x)$$ where each $f^a_b,g^a_b$
has degree $b$ and the following relations, equivalent to flatness,
hold $$yf=vg$$ $$xg=uf$$\end{prop}\proof See \cite {R2}\endproof\par
Concretely, the above relation mean\nl - the coefficients of
$f^1_{m-i}(x), g^1_{i-1}(y)$ are free parameters (no relations);\nl
- the relation $uv=t$ holds;\nl - $f^2_{i-2}, g^2_{m-i-1}$ are
determined by the other data.\ss A similar and simpler story holds
at the principal ideals $I^m_i(a)$.\par We conclude \item - $X\sbr
m._B$ is smooth;\item - its fibre at $t=0$, i.e. $\Hilb_m(X_0)$ has,
along $\Hilb^0_m(R)$, $(m+1)$ smooth components crossing normally,
$D_0,...,D_m.$ \vskip 3cm\begin{center}fig. 2\end{center}\vskip 1cm
In fact, if $X_0=X'_0\cup X''_0$ then
$$D_i\sim_{bir}(X'_0)^{m-i}\times(X''_0)^i.$$ The next question is:
how to glue
together the various local deformations ?\ss
{\underline{\bf{Construction}}} Let $C_1,...,C_{m-1}$ be copies of
$\P^1$, with homogenous coordinates $u_i,v_i$ on the $i$-th copy.
Let $\tilde{C}\subset C_1\times...\times C_{m-1}\times B$ be the
subscheme defined by\ms
$$v_1u_2=tu_1v_2,...,
v_{m-2}u_{m-1}=tu_{m-2}v_{m-1}$$\ms Fibre of $\tilde{C}$ over
$0\in B$ is
$$\tilde{C}_0=\bigcup\limits_{i=1}^m\tilde{C}_i,$$
where $$\tilde{C}_i= [1,0]\times...\times[1,0]\times
C_i\times[0,1]\times...\times[0,1]$$ In a neighborhood of
$\tilde{C}_0$, $\tilde{C}$ is smooth and $\tilde{C}_0$ is its
unique singular fibre over $B.$ We may embed $\tilde{C}$ in
$\P^{m-1}\times B$ via $$Z_i=u_1\cdots u_{i-1}v_{i}\cdots v_{m-1},
i=1,...,m.$$ These satisfy
$$Z_iZ_j=tZ_{i+1}Z_{j-1}, i<j-1 $$ so embed $\tilde{C}$ as a
family of rational normal curves $\tilde{C}_t\subset\P^{m-1},
t\neq 0$ specializing to a connected $(m-1)$-chain
 of lines.\par
 Next consider  $\A^{2m}$ with coordinates $a_0,...,a_{m-1},
d_0,...,d_{m-1}$ \nl Let $\tilde{H}\subset\tilde{C}\times\A^{2m}$
be defined by\ss  $\begin{matrix} &a_0u_1=tv_1,&\\
a_1u_1=d_{m-1}v_{1},&\ldots
&,a_{m-1}u_{m-1}=d_1v_{m-1}\\&d_0v_{m-1}=tu_{m-1}&
\end{matrix}$

 Fibres of
$\tilde{H}$ over $\A^{2m}$ are: a point (generically), or a chain
of $i\leq m-1$ rational curves; all values $i=1,...,m-1$ occur.
Consider the subscheme of $Y=\tilde{H}\times_{B}U$ defined by
\begin{eqnarray*}
F_0:=x^m+a_{m-1}x^{m-1}+...+a_1x+a_0~~~~~~~~~~~~~~\\
F_1:=u_1x^{m-1}+u_1a_{m-1}x^{m-2}+...+u_1a_2x+u_1a_1\\
+v_1y~~~~~~~~~~~~~~~~~~~~~~~~~~~~~~~~~~~~~~~~\\
 ...~~~~~~~~~~~~~~~~~~~~~~\\
F_i:=u_ix^{m-i}+u_ia_{m-1}x^{m-i-1}+...+u_ia_{i+1}x+u_ia_i\\ +
v_id_{m-i+1}y+...+v_id_{m-1}y^{i-1}+ v_iy^i ~~~~~~~~~~~~\\
 ...~~~~~~~~~~~~~~~~~~~~~~~\\
F_m:=d_0+d_1y+...+d_{m-1}y^{m-1}+y^m~~~~~~~~~~~\end{eqnarray*}

\ss The following is proven in \cite{R3}
\begin{thm} (i) $\tilde{H}$ is smooth and irreducible.\nl (ii)
The ideal sheaf $\I$ generated by $F_0,...,F_m$ defines a subscheme
of $\tilde{H}\times_BU$ that is flat of length $m$ over $\tilde{H}$
\nl (iii)The classifying map $$\Phi=\Phi_\I:
\tilde{H}\to\Hilb_m(U/B)$$ is an isomorphism.\end{thm} The proof
shows furthermore that $\tilde{H}$ is covered by opens
$$U_i=\{Z_i\neq 0\}, i=1,...,m$$\vskip
3cm\begin{center}fig.3\end{center}\vskip 1cm On $U_i,$ we have
$$F_j=u_jx^{i-j-1}F_{i-1}, j<i-1$$ $$F_j=v_jy^{j-i}F_i, j>i$$ hence
$ F_{i-1}, F_i$ generate $\I$ on $U_i$ (they yield the $f,g$ in the
universal deformation of Proposition 2 above).\par Also,
$a_i=(-1)^i\sigma^x_{m-i}$ are the elementary symmetric functions in
the roots of $F_0$, and ditto for $d_i, \sigma^y_{m-i}, F_m.$ So the
projection $\tilde{H}\to\A^{2m}_B$ factors through the cycle map
$$\begin{matrix}\ \ \tilde{H}&&\\\frak c\downarrow&\searrow&
\\ \ \ \ \ X\sfl
m._B&\stackrel{\sigma}{\to}&\A^{2m}_B\end{matrix}$$
$$\sigma=(\sigma^x_1,...,\sigma^x_m,\sigma^y_1,...,\sigma^y_m)$$
(one can show $\sigma$ is embedding).To prove the Main Theorem, we
must show: $\frak c$ is the blow-up of $D^m$\nl It is convenient to
pass to an 'ordered' model, defined by the following Cartesian
diagram:
$$\begin{matrix}X^{\lceil m\rceil}_B&\to &X\sbr m._B
\\\downarrow&&\downarrow\\X^m_B&\to&X\sfl
m._B\end{matrix}$$ In this diagram, the right vertical arrow is the
cycle map, the bottom horizontal arrow is the natural map between
the Cartesian and symmetric products, and the other arrows are
defined by the fibre product construction.
 Recall
the description of the blowup of an ideal $I$ as subscheme of
$\P(I)$. Let us rewrite the defining local equations for $X\sbr
m._B$ in terms of the homogeneous coordinates $Z_i$ on $\P^{m-1}$:
they are
\begin{flushleft}{\emph{linear}}:
$$\sigma_{m-j}^yZ_i=t^{m-j-i}\sigma_j^xZ_{i+1},\
 i=1,...,m-1,j=0,...,m-1; $$
$$\sigma_{m-j}^xZ_i=t^{m-j-i}\sigma_j^yZ_{i-1},\
i=2,...,m, j=0,...,m-1. $$ {\emph{quadratic}}:
$$Z_iZ_j=tZ_{i+1}Z_{j-1}, i<j-1 $$\end{flushleft}
Our task at this point is to 'reverse engineer' an ideal whose
generators $G_1,...,G_m$ satisfy (precisely) these relations.
Actually, the choice of $G_1$ determines $G_2,...,G_m$ via the
linear relations, though a priori, $G_2,...,G_m$ are only
\emph{rational} functions. Now recall that $Z_1$ generates $\O(1)$
over the open $U_1$ which meets the special fibre $t=0$ in the
locus of $m$-tuples entirely on $x$-axis. On that locus, an
equation for the discriminant is given by the Van der Monde
determinant:
$$v^m_x=\det(V^m_x),$$
$$V^m_x=\left [ \begin{matrix}
 1&\ldots&1\\x_1&\ldots&x_m\\\vdots&&\vdots\\
x_1^{m-1}&\ldots&x_m^{m-1}
\end{matrix} \right ].$$
Thus motivated, set $$G_1=v^m_x.$$ This forces
$$G_i=\frac{(\sigma_m^y)^{i-1}}{t^{(i-1)(m-i/2)}}v^m_x=
\frac{(\sigma_m^y)^{i-1}}{t^{(i-1)(m-i/2)}}G_1,\ \ i=2,...,m.$$ If
the construction is to make sense, these better be regular. In
fact,
$$G_i=\pm\det(V^m_i),$$
$$V^m_i=\left[\begin{matrix}
 1&\ldots&1\\x_1&\ldots&x_m\\\vdots&&\vdots\\
x_1^{m-i}&\ldots&x_m^{m-i}\\y_1&
\ldots&y_m\\\vdots&&\vdots\\
y_1^{i-1}&\ldots&y_m^{i-1}
\end{matrix}\right]  $$
(we call this the 'Mixed' Van der Monde matrix). The $G_i$ satisfy
same relations as the $Z_i$, so we can map isomorphically
$$\O(1)\to J=\Ideal(G_1,...,G_m)$$
$$Z_i\mapsto G_i.$$Then  $J$ is an invertible ideal defining a Cartier
divisor $\Gamma$. The Main Theorem's assertion that $\frak c$ is
the blowup of $D^m$ means
$$J=\frak c^*(I_{D^m})$$ i.e. $$\Gamma=\frak c^*(D^m)$$
Containment $\supseteq$ is clear. Equality is clear off the
special fibre $t=0$. Now this special fibre is  sum of components
$$\Theta_I=\text{Zeros}(x_i,i\not\in I, y_i, i\in I),
 I\subseteq\{1,...,m\}.$$
 Set
 $$\Theta_i=\bigcup\limits_{|I|=i}\Theta_I.$$
 Note that the open set $U_i$ meets only $\Theta_i, \Theta_{i-1}$.
 One can check that the vanishing order of
 $G_j$ on any $\Theta_I, |I|=k,$  is $$\ord_{\Theta_I}(G_i)=(k-i)^2+(k-i)$$ $$=0\ \ \
 {\text{if}}\ \ \  k=i, i-1$$ So $\Zeros(G_i)=\frak c^*(D^m)$ on
 $U_i$, i.e. $$\Gamma|_{U_i}=\frak c^*(D^m)|_{U_i},\ \ \forall i$$
 $$\therefore\ \ \  \Gamma=\frak c^*(D^m)$$ This concludes the proof of the
 Main Theorem.\par One point of interest is the interpretation
 of mixed Van der Monde matrices, whose determinants played a large role
 in the proof :\ The universal subscheme
 $$\Xi=\Zeros(\I)\subset X^{\lceil m\rceil}_B\times X$$ contains sections
 $$\Psi_i={\mathrm{graph}}(p_i:X^{\lceil m\rceil}_B\to X)$$ The universal
 quotient $$\frak Q_m=p_{X^{\lceil m\rceil}_B*}(\O/\I)$$ maps to $\O_{X\sfl
 m._B}$ via restriction on $\Psi_i$. Assembling together, get map
 $$V:\frak Q_m\to m\O_{X^{\lceil m\rceil}_B}.$$ Then $V^m_i=$ is just the
 matrix of $V$
 with respect to the basis $1,x,...,x^{m-i}, y,...,y^{i-1}$
 of $\frak Q_m$ on
 $U_i.$\par
 A somewhat mysterious point that comes up in the above proof is:
 as the $Z_i$ are interpreted as the equations of the discriminant,
 what,if any,
is the interpretation of $u_i, v_i$ ?
\section{Applications}\subsection*{Canonical bundle}
 A first application is a formula for the canonical bundle of $X\sbr m._B$.
  For
 any class $\alpha$ on $X$, denote $$\alpha\sbr m.=q_*p^*(\alpha)$$
 where $\Omega\subset X\sbr m._B\times_BX$ is universal subscheme
 and $p:\Omega\to X, q:\Omega\to X\sbr m._B$ are natural
 maps. Here $q_*$ denotes the cohomological direct image, sometimes
 called the norm or denoted $q_!$, \emph{not} the sheaf-theoretic
 direct image.\par Another way to construct $\alpha\sbr m.$
 is as follows. First note the natural isomorphism over $\Q$
\begin{equation*}H^*(X\spr m.)\simeq\Sym^m(H^*(X))\end{equation*}
This yields a class $\alpha^m\in H^*(X\spr m.)$, and $\alpha\sbr
m.$ is the image of the latter via the composite
 \begin{equation}
 H^*(X\spr m.)\to H^*(X\spr
 m._B)\stackrel{\frak c_m^*}{\to} X\sbr m._B\end{equation}
 Also set $$\O_{X\sbr m._B}(1)=\O(-\Gamma)$$ (the canonical $\O(1)$ as
 blowup, via $\Proj$).
 \begin{cor}\label{khilb}
 {$K_{{X\sbr m._B}/B}=(K_{X/B})\sbr m. \otimes\O_{X\sbr m._B}(1)$}
 \end{cor}
 \proof It suffices to note that  both sides
 agree off the exceptional locus of $\frak c_m$.\endproof
 In particular, $K_{{X\sbr m._B}/B}.C^m_i=+1$, so $\frak c_m$
 'looks like' a flip. The following example partially confirms
 this.
 \subsection*{Example: $m=2$.} We have a diagram
  $$\begin{matrix}
 X^{\lceil 2\rceil}_B&\to &X^{[2]}_B\\
 \frak c'_2\downarrow&&\downarrow\frak c_2\\
 X^2_B&\to&X\sfl 2._B
 \end{matrix}$$ with horizontal maps of degree 2. Local equations
  for $X^2_B$ are:$$ x_1y_1=x_2y_2=t$$ (so this is a 3-fold
 ODP);\nl for $X^{\lceil 2\rceil}_B$: $$x_1u=y_2v$$$$x_2u=y_1v$$
 so $\frak c'_2$ is a small resolution of the ODP, known as a flopping
 contraction; it can be flopped to
 yield $X^{**}$ smooth that is the source of the 'opposite' flopping contraction.\nl
 Equations for $X\sfl 2._B$ are: $$\sigma_2^y\sigma_1^x=t\sigma_1^y$$
 $$\sigma_2^x\sigma_1^y=t\sigma_1^x$$
 $$\sigma_2^x\sigma_2^y=t^2$$ This is a cone over a cubic scroll
 in $\P^4.$\  $X^{[2]}_B$ is small resolution of the cone, with
 exceptional locus $C^2_1=\P^1.$ A well-known procedure, due to
  Francia, yields a flip, called Francia's flip, of $\frak c_2$: blow up
 $C^2_1$ in $X^{[2]}_B$ (which is the same as blowing up the
 vertex of the cone);
 the exceptional divisor is a scroll of type $ F_1$; then
 blow up the negative curve of $F_1$ to get a new exceptional surface
 of type $F_0$; then
 blow down $F_0$ in the other direction to  $C^*=\P^1$ so the $F_1$
 becomes a $\P^2$; then finally blow down $\P^2$ to a
 (singular) point on a new
 3-fold $X^*$,
which  is 2:1 covered by $X^{**}.$\par This situation is
intriguing in view of recent work of Bridgeland \cite{B} and
Abramovich and Chen \cite{AC} which shows that the flop $X^{**}$
and the flip $X^*$ can be interpreted as moduli spaces of certain
'1-point perverse sheaves' on $X^{\lceil 2\rceil}_B$ and $X\sbr
2._B$, respectively. This raises the question of finding a natural
interpretation of
 $X^*, X^{**}$ and their higher-order analogues, if they exist,
  in terms of our family of curves $X/B$.
\ss
  \subsection*{Euler number} As an application of our study of $\frak c_m$,
 we can
 compute (topological) Euler number $e({X\sbr m._B})=c_{m+1}(T_{X\sbr
 m._B})$, at least for case of $\leq 1$ node in any fibre:
 \begin{cor} If $X/B$ has $\sigma$ singular fibres and each has
 precisely 1 node, then the topological Euler number of $X\sbr
 m._B$ is given by
 \begin{eqnarray}\label{euler}e({X\sbr m._B})
 =(-1)^m\binom{2g-2}{m}(2-2g(B)) +\sigma\binom{m-2g+2}{m-1}
 \end{eqnarray}
 \end{cor}
\proof Let $$(X_i, p_i, X_{i,0}=X_i\setminus p_i),
i=1,...,\sigma$$ be the singular fibres with their respective
unique singular point and smooth part, and
$$X_0=X\setminus (X_1\cup\ldots\cup X_\sigma), B_0=\pi(X_0).$$
Then $X\sfl m._B$ admits a (locally closed) stratification with
big stratum $$(X_0)\sfl m._{B_0}$$ and other strata
$$\Sigma_{i,j}=ip_j+(X_{j,0})\sfl {m-i}. , i=0,...,m,
j=1,...,\sigma.$$ The fibre of $\frak c_m$ over each of these strata
is, respectively, a point over the big stratum, and over the
$\Sigma_{i,j},$ a point for $i=0,1$, a chain of $(i-1)$ $\P^1$\ s
for $i=2,...,m.$ Since the Euler number is multiplicative in
fibrations and additive over strata, we get
\begin{eqnarray*}e({X\sbr m._{B}})=e((X_0)\sfl m._{B_0})+
\sum e((X_{j,0})\sfl m.)\\
+\sum\limits_{i>0}ie((X_{j,0})\sfl {m-i}.)
\end{eqnarray*}
Now MacDonald's formula \cite{Mac} says that for any $X$, the
Euler number of its $m$th symmetric product is given by
$$e(X\sfl m.)=(-1)^m\binom{-e(X)}{m}.$$ Plugging this into the
above and using multiplicativity for the fibration $(X_0)\sfl
m._{B_0}$ over $B_0$ yields
\begin{equation}\label{euler2}
e({X\sbr
m._B})=(-1)^m\binom{2g-2}{m}(2-2g(B))~~~~~~~~~~~~~~~~~~~~~~~~~~~~
 ~~~~~~~~~~~~~~~~~~~~~~~+\sigma\sum\limits_{k=0}^{m-1}
(-1)^k(m-k)\binom{2g-2}{k}\end{equation} Now, as pointed out by
L.C. Wang, (\ref{euler}) follows from (\ref{euler2}) by the
elementary formula
$$\sum\limits_{k=0}^b (-1)^k\binom{a}{k}=(-1)^b\binom{a-1}{b}$$
which in turn is an easy consequence of Pascal's relation\nl
$\binom{a}{k}=\binom{a-1}{k}+\binom{a-1}{k-1}.$\endproof

\begin{rem}\rm  Suppose our family $X/B$ is a blowup
$$\beta:X\to Y$$
of  a \emph{smooth} $\P^1$ bundle; equivalently, each singular fibre
of $X/B$ has consists of two $\P^1$ components. Then there is
another way to construct $X\sbr m._B$ and obtain formula
(\ref{euler2}) above, as follows.
 Note that the natural map
 $$\eta:Y\sbr m._B=Y\sfl m._B\to B$$ is a $\P^m$-bundle. Blow
 up a $\P^{m-1}$ in each  fibre of $\eta$ over a singular value
 of $\pi$, giving rise to exceptional divisors $E_{1,i},
 i=1,...,\sigma$;
 then blow up a $\P^{m-2}$ in general position in each exceptional
 divisor $E_{1,i},$ giving rise to new exceptional divisors $E_{2,i},$
 etc. Finally, blow up general point on each exceptional
 divisor $E_{m-1,i}$. This yields $X\sbr m._B$.
 In these blowups, the change in Euler number is easy to
 analyze, yielding (\ref{euler}).\qed
\end{rem}

 \subsection*{Further developments (under construction)}
\vskip 1cm
  We mention
 some natural questions and possible extensions.

 $\bullet$ What is the total Chern class $c(T_{X\sbr m._B}) $?\par
 $\bullet$ Develop intersection calculus for
  diagonal loci of all codimensions in
 $X\sbr m._B$, i.e. degeneracy loci $$\Gamma^m_r=\rk(V^m_i)\leq m+1-r$$
 (locus where $r$ points come together)\nl More generally, loci
 $\Gamma^m_{(m.)}$, $m_1+...+m_k=m$,
 $$\Gamma^m_{(m.)}=\{z:\frak c_m(z)=\sum m_ip_i\}.$$ In particular,
 the small diagonal
 $$\Gamma^m_{(m)}=\text{locus of length-$m$ schemes supported at 1
 point}$$ which coincides with the blowup of $X$, locally at each fibre
 node, in a punctual subscheme of type $$(x^{\binom{m}{2}},...,
 x^{\binom{m-i}{2}}y^{\binom{i}{2}},...,y^{\binom{m}{2}})$$
A potential application of this calculus is is to enumerative
geometry (multiple points, multisecant
 spaces, special divisors on stable curves...)\nl
A \emph{Sample corollary} which however can also be derived by
other means) is the following \textbf{ relative triple point
 formula:} for a map $f:X\to\P^2$, the number of relative triple
 points is
  \begin{eqnarray*} N_{3,X}(f)=
 (\frac{(d-2)(d-4)}{2}+g-1)L^2+(3-\frac{d}{2})\omega
 L+2\omega^2-4\sigma\end{eqnarray*}
 \nl where $L=f^*\O(1), L^2=\deg(f),$
 $ d=\deg($fibre), $ g=$genus(fibre).\nl See\cite{R3} for some progress on this.
 \par
 $\bullet$ If $X/B$ is of
 \emph{compact type} (assume for simplicity there exists a section),
 we have an Abel-Jacobi morphism to the Jacobian: $$X^{[m]}_B\to
 J(X/B)$$
 Fibres give a notion of 'generalized linear system' on reducible
 fibres. How is this related to other approaches to such notions in the
 literature ?
 \section{Chern classes of tautological bundles}
 In \cite{R} we gave a simple formula for the Chern classes of the
 tautological bundles $\lambda_m(L)$, where $L$ is a vector bundle
 on $X$. Here $X$ need not be a surface; we just need
 a family of nodal curves $X/B$. More precisely,
 we gave in \cite{R} a formula for the pullback of
 $\lambda_m(L)$ on the (full) flag relative Hilbert scheme,
 denoted $W^m(X/B).$ The formula is simple and involves only
 divisor classes plus classes coming from $X$,
 but has the disadvantage that these classes,
 unlike $\lambda_m(L)$ itself, do not descend to the Hilbert
 scheme $X\sbr m._B$. Though it is, broadly speaking,
  obvious that a formula on
 $X\sbr m._B$ can be derived from the one on $W^m(X/B),$ it is still of
 some interest, in view of possible applications, to work this out.
 It turns out that for $X$ a surface,
  a formula for the Chern classes of tautological
 bundles was already derived, in the context of the (absolute)
 Hilbert scheme $X\sbr m.$, by Lehn \cite{L}, using the Fock
 space formalism introduced earlier by Nakajima \cite{N, EG}.
 Since our tautological bundles $\lambda_m(L)$ are pullbacks of
 the analogous bundles on $X\sbr m.$ via the natural inclusion
 $$X\sbr m._B\subset X\sbr m.,$$ Lehn's formula yields an
 analogous one on $X\sbr m._B.$ Our purpose here, then, is to
 verify that when $X$ is a surface,
  the push-down from $W^m(X/B)$ to $X\sbr m._B$ of the
 formula of \cite{R} coincides with the restriction of Lehn's
 formula, at least when $L$ is a line bundle. Thus, we have
 compatibility in the natural diagram
 $$\begin{matrix}
 W^m(X/B)&&\\ w_m\downarrow&&\\ X\sbr m._B&\to&X\sbr m.\end{matrix}
 $$
 \par We begin with some formalism. First, we have the
 operation of \emph{exterior multiplication} $\star$ of cohomology
 classes on various $X\sbr m._B$, defined as follows. Let
 $$Z_{m,n}\subset X\sbr m._B\times_B X\sbr n._B\times_B X\sbr
 m+n._B$$ be the closure of the locus $$(z_m, z_n, z_m\coprod
 z_n)$$ (where $z_m, z_n$ are disjoint), and let
 $$p: Z_{m,n}\to X\sbr m._B\times_B X\sbr n._B,
 q:Z_{m,n}\to X\sbr m+n._B$$ be the projections,
 both generically finite. For $\alpha\in
 H^r(X\sbr m._B), \beta\in H^s(X\sbr n._B),$
 identifying homology and cohomology, set
 $$\alpha\star\beta=q_*p^*(\alpha\times\beta)\in H^{r+s}(X\sbr
 m+n._B).$$ This operation is obviously associative and
 commutative on even (in particular, algebraic) classes.
 In particular, taking $\beta=1,$ we get a natural way of mapping
$H^r(X\sbr m._B)$ to
 $H^r(X\sbr m+n._B)$ for each $n\geq 0$.
 \par Next, consider the small diagonal
 $$\Gamma^m_{(m)}\stackrel{i_m}{\hookrightarrow} X\sbr m._B.$$ The
 restriction of the cycle map yields a birational morphism
 $$\beta_m:\Gamma^m_{(m)}\to X.$$ For any $\alpha\in H^r(X),$ we set
 $$q_m[\alpha]=i_{m*}(\beta_m^*(\alpha))\in H^{r+2m-2}(X\sbr m.),$$
 Via $\star$ multiplication, $q_m[\alpha]$ may be viewed as
 with an operator on $$\bigoplus\limits_{s,n=0}^{\infty}H^s(X\sbr m._B)$$
 which has operator bidegree $(r+2m-2, m).$ This is
 known as Nakajima's creation operator (cf. \cite{ N, EG}).\par
  Lehn's formula is as follows
 \begin{thm*}(Lehn \cite{L}) For a line bundle $L$,
 the total Chern class of $\lambda_m(L)$ is the part
 in bidegrees $(*,m)$ of $$\exp(\sum\limits_{n=1}^\infty
 \frac{(-1)^{n-1}}{n}q_n[c(L)]).$$

\end{thm*}

Now our formula is the following
\begin{thm}\label{tauto} For a line bundle $L$, we have\nl
\lefteqn{ c(\lambda_m(L))=}
\begin{equation*}\sum\limits_{\begin{tiny}
\begin{matrix}I=(1\leq i_1<...<i_k)\\
|I|\leq m\end{matrix}\end{tiny}}
(-1)^{|I|-k}\frac{(i_1-1)!...(i_k-1)!}{|I|!(m-|I|)!}
q_{i_1}[c(L)]...q_{i_k}[c(L)].\end{equation*}
\end{thm}

It is elementary to derive Theorem \ref{tauto} from Lehn's theorem
(whose proof is rather long). Our purpose here, however, is to
derive Theorem 2 from a result in \cite{R}, as follows. Let
$$w_m:W=W^m(X/B)\to X\sbr m._B$$ be the natural morphism from the
flag Hilbert scheme to the ordinary one, let $$p_i:W\to X$$ be the
$i$th projection, mapping a filtered scheme $z_1<...<z_m$ to the
support of $z_i/z_{i-1},$ and let
$$\Delta_{ij}\subset W, i<j$$
denote the (reduced) locus where the $p_i$ and $p_j$ coincide; also
set, for any class $c\in H^*(X),$  $$c_i=p_i^*(c).$$ It is shown in
\cite {R} that each sum $\sum\limits_{i=1}^{j-1}\Delta_{ij}$ is a
Cartier divisor (even though $W$ is in general singular and each
summand individually is not Cartier). It is also shown there that
the following result holds (for a line bundle $L$):
\begin{eqnarray}\label{geo}
c(w^*\lambda_m(L))=\prod\limits_{j=1}^m(1+L_j-
\sum\limits_{i=1}^{j-1}\Delta_{ij})
\end{eqnarray}
Deriving Theorem \ref{tauto} from (\ref{geo}) is a matter of
expanding the product as a sum of monomials, applying $w_*$ and
dividing by $m!=\deg(w).$ In doing so, it is useful to observe the
following. Let's call a \emph{connected} monomial  on an index set
$I$ one which, after a permutation, can be written in the form
$$q_I[c]=c_{i_1}\Delta_{i_1i_2}\Delta_{i_2i_3}...\Delta_{i_{k-1}i_k},
I=(i_1<...<i_k)$$ where $c$ is either 1 or $[L]$. The intersection
implicit in the above product is transverse, hence well-defined even
though the divisors are not Cartier.
 It is easy to
see by induction that there are $(k-1)!$ unordered monomials in
the expansion of (\ref{geo}) yielding the same $q_I[c]$. Moreover
it is clear that
$$w_*(q_I[c])=q_{|I|}[c].$$ Now we note that each monomial
appearing in the expansion of (\ref{geo}) may be decomposed
uniquely as a product of connected monomials on pairwise disjoint
index sets (its 'connected components'), yielding a term
$$(-1)^{\sum\limits_{j=1}^k(|I_j|-1)}q_{I_1}[c_1]...q_{I_k}[c_k],$$
each $c_j\in\{1,[L]\}$ which, for fixed $I=I_1\coprod ...\coprod
I_k$, appears $(|I_1|-1)!...(|I_k|-1)!$ times. Applying $w_*$, we
get, for each choice of $I\subseteq \{1,...,m\}$ and $c_j$, a term
in $w_*$ applied to (\ref{geo}):
$$(i_1-1)!...(i_k-1)!(-1)^{i-k}q_{i_1}[c_1]...q_{i_k}[c_k],$$
$i_j=|I_j|, i=\sum\limits_ji_j.$ Then multiplying by
$\binom{m}{i}$ for the choice of subset $I$ with $|I|=i$, and
dividing by $m!$ yields the result.\qed

\end{document}